\documentclass[11pt]{article}
\usepackage{amsmath}
\usepackage{amsmath,amssymb,latexsym,color}
\usepackage[mathscr]{eucal}
\usepackage{amssymb,amsmath}

\newtheorem{thm}{Theorem}[section]

\newtheorem{lemma}[thm]{Lemma}

\newtheorem{example}{Example}[section]
\newtheorem{defin}{Definition}[section]

\newtheorem{remark}{Remark}[section]

\newcommand{\proof}{{\it Proof.\quad}}
\newcommand{\qed}{\hfill\Box\medskip}

\usepackage{CJK}
\begin{document}
\begin{CJK*}{GBK}{song}

\newcommand{\be}{\begin{equation}\label}
\newcommand{\ee}{\end{equation}}
\newcommand{\bea}{\begin{eqnarray}\label}
\newcommand{\eea}{\end{eqnarray}}

\title{\bf The smallest one-realization of a given set \uppercase\expandafter{\romannumeral3}}

\author{
Kefeng Diao$^{\rm a}$\quad Ping Zhao$^{\rm a}$\quad    Kaishun
Wang$^{\rm b}$\thanks{Corresponding
author: wangks@bnu.edu.cn}\\
{\footnotesize a. \em School of Science, Linyi  University,
Linyi, Shandong, 276005, China }\\
{\footnotesize b. \em  Sch. Math. Sci. {\rm \&} Lab. Math. Com.
Sys., Beijing Normal University, Beijing 100875,  China} }

\date{}
 \maketitle

\begin{abstract}

In [The smallest one-realization of a given set,  Electronic J.
Combin. 19 (2012), $\sharp$P19],   we determined the minimum number
 of vertices of one-realizations of a given finite set $S$, and  constructed    the
corresponding mixed hypergraphs. In this paper, by finding some of
their spanning sub-hypergraphs, we determine the minimum number of
${\cal D}$-deges (resp. ${\cal C}$-edges) of one-realizations of
$S$. As a result, we partially solve an open problem proposed by
Tuza and Voloshin in 2008.

\medskip
\noindent {\em Key words:} hypergraph coloring; mixed hypergraph;
feasible set; one-realization.
\end{abstract}

\section{Introduction}

 A {\em mixed hypergraph } on a finite set $X$ is a triple ${\cal
H}=(X, {\cal C}, {\cal D})$, where ${\cal C}$ and ${\cal D}$ are
families of subsets of $X$, called the {\em ${\cal C}$-edges} and
{\em ${\cal D}$-edges}, respectively. A \emph{bi-edge} is an edge
which is both a ${\cal C}$-edge and a ${\cal D}$-edge. If ${\cal
C}={\cal D}$, ${\cal H}$ is a {\em bi-hypergraph}. If each edge has
$r$ vertices, then ${\cal H}$ is {\em $r$-uniform}. When $\cal
C=\emptyset$ (resp. $\cal D=\emptyset$), $\cal H$ is a ${\cal
D}$-{\em hypergraph} (resp. ${\cal C}$-{\em hypergraph}).
 A sub-hypergraph ${\cal H}'=(X', {\cal C}', {\cal D}')$ of  ${\cal H}=(X, {\cal C}, {\cal D})$  is a {\em spanning
sub-hypergraph}  if $X'=X$, and
 ${\cal H}'$ is   a {\em derived sub-hypergraph}
  of  ${\cal H}$ on $X'$, denoted by ${\cal H}[X']$, when ${\cal C}'=\{C\in {\cal C}| C\subseteq X'\}$
 and ${\cal D}'=\{D\in {\cal D}| D\subseteq X'\}$.


   A {\em proper $k$-coloring} of
${\cal H}$ is a mapping from $X$ into a set of $k$ {\em colors} so
that each ${\cal C}$-edge has two vertices with a {\em Common} color
and each ${\cal D}$-edge has two vertices with {\em Distinct}
colors. A {\em strict $k$-coloring} is a proper $k$-coloring using
all of the $k$ colors, and a mixed hypergraph is {\em $k$-colorable}
if it has a strict $k$-coloring.  A coloring of   ${\cal H}$ may  be
viewed as a  partition  of its vertex set, where the \emph{color
classes}    are the sets of vertices assigned to the same color, so
a strict $n$-coloring $c=\{C_1,C_2,\ldots,C_n\}$ of  ${\cal H}$
means that  $C_1,C_2,\ldots, C_n$ are  the $n$ color classes under
$c$. The set of all the  values $k$ such that ${\cal H}$ has a
strict $k$-coloring is called the {\em feasible set}
  of ${\cal H}$, denoted by ${\cal F}({\cal H})$. For each $k\in {\cal F}({\cal H})$, let $r_k$ denote
the number of {\em partitions} of the vertex set. For  a set $S$ of
positive integers, we say that a mixed hypergraph ${\cal H}$ is a
{\it realization} of $S$ if ${\cal F}({\cal H})=S$. A mixed
hypergraph ${\cal H}$ is a {\it one-realization} of $S$ if it is a
realization of $S$ and $r_k=1$ for each $k\in S$.  When one
considers the colorings of a mixed hypergraph, it suffices to assume
that each ${\cal C}$-edge   has at least three vertices.  The study
of   the colorings of mixed hypergraphs has made a lot of progress
since its inception \cite{Voloshin1}. For more information, we would
like refer readers to \cite{Kobler, Tuza, Voloshin2, Voloshin3}.

K$\ddot{u}$ndgen et al. \cite{Kundgen} initiated  a systematic study
of planar mixed hypergraphs, and found a one-realization of
$\{2,4\}$ on 6 vertices for planar hypergraphs. Bujt$\acute{\rm a}$s
and Tuza \cite{Bujtas} gave a necessary and sufficient condition for
a set to be the feasible set of an $r$-uniform mixed hypergraph.
 In \cite{zdw},   we introduced a new construction of hypergraphs, and characterized the feasible set and
chromatic spectrum of a $3$-uniform bi-hypergraph. Jiang et al.
\cite{Jiang} proved that a set $S$ of positive integers is a
feasible set of a mixed hypergraph if and only if $1\notin S$ or $S$
is an interval, and determined the minimum number of vertices of
realizations of $\{s,t\}$ with $2\leq s\leq t-2$. Kr$\acute{\rm
a}$l \cite{Kral} initiated the study of the minimum number of
vertices of one-realizations of $S$, and obtained an upper bound. In
\cite{zdw1}, we determined this minimum number. We further determined the minimum number of vertices of
3-uniform bi-hypergraphs which are one-realizations of $S$ in
\cite{zdw2}.

Note that the feasible set  of any  $r$-uniform ${\cal
C}$-hypergraph contains $\{1,\ldots,r-1\}$ as a subset. On $n$
vertices, the minimum number of $r$-element ${\cal C}$-edges to
generate this smallest possible feasible set is $\lceil
n(n-2)/3\rceil$ in the particular case of $r=3$, but only some lower
and upper estimates of the order $\Theta(n^{r-1})$ are known if
$r\geq 4$ (see \cite{Diao1, Diao2}). So Tuza and Voloshin
\cite{Tuza} proposed the following problem:

\textbf{Problem.} Given a finite set $S$ of positive integers,
determine or estimate the  minimum numbers of (${\cal
C}$-, ${\cal D}$-, bi-) edges in a mixed (bi-) hypergraph whose
feasible set is $S$.

In order to solve this problem for one-realizations of a given set,
we need to deleting some $\cal D$-edges or $\cal C$-edges from the
mixed hypergraphs constructed in \cite{zdw1}. Now we introduce the
mixed hypergraphs.

In the rest we always assume that $S=\{n_1,n_2,\ldots,n_s\}$ is a
set of integers with $2\leq n_s<\cdots <n_2<n_1$ and  $[n]$ is the
set $\{1,2,\ldots,n\}$.

\medskip
\noindent {\bf Construction} (\cite{zdw1}) Let
 \begin{eqnarray*}X_{n_1,\ldots,n_s}&=&\{(n_1,n_2,\ldots,n_s)\}\cup\{(\underbrace{i,i,\ldots,i}_s)|~i\in[n_s-1]\}\\
&\cup&\bigcup_{t=2}^{s}\bigcup_{j=n_t}^{n_{t-1}-1}
\{(\underbrace{j,\ldots,j}_{t-1},n_t,n_{t+1},\ldots,n_s),
(\underbrace{j,\ldots,j}_{t-1},\underbrace{1,\ldots,1}_{s-t+1})\},\\
 {\cal
D}_{n_1,\ldots,n_s}&=&\{\{(x_1,x_2,\ldots,x_s),(y_1,y_2,\ldots,y_s)\}|x_i\neq y_i, i\in [s]\},\\
 {\cal
C}_{n_1,\ldots,n_s}&=&\{\{(x_1,\ldots,x_s),(y_1,\ldots,y_s),(z_1,\ldots,z_s)\}|~
|\{x_j, y_j, z_j\}|=2, j\in [s]\}.
\end{eqnarray*}
By \cite[Theorems 1.1, 2.4 and 2.5]{zdw1} the  hypergraph ${\cal
H}_{n_1,\ldots,n_s}= (X_{n_1,\ldots,n_s},{\cal C}_{n_1,\ldots,n_s},
{\cal D}_{n_1,\ldots,n_s})$ (resp.  ${\cal G}_{n_1,\ldots,n_s}={\cal
H}_{n_1,\ldots,n_s}[X_{n_1,\ldots,n_s}\setminus
\{(n_2,1,\ldots,1)\}]$) is a smallest one-realization of $S$ when
$n_1-1\not\in S$ (resp. $n_1-1\in S$).

 By finding  some spanning
sub-hypergraphs of ${\cal H}_{n_1,\ldots,n_s}$ and  ${\cal
G}_{n_1,\ldots,n_s}$, we determine the minimum number of ${\cal
D}$-deges (resp. ${\cal C}$-edges) of one-realizations of $S$ as
follows.

\begin{thm}\label{thm1} Let $\delta_{{\cal D}}(S)$ denote the minimum number of ${\cal D}$-edges of one-realizations of $S$. Then
  $$\delta_{{\cal D}}(S)=\left \{
\begin{array}{ll}\frac{n_1(n_1-1)}{2}, & \mbox{if~ $n_1-1\notin S,$}\\
\frac{n_1(n_1-1)}{2}-1,&\mbox{if ~$n_1-1\in S.$}\\
\end{array}
\right. $$
\end{thm}

\begin{thm}\label{thm2} Let $\delta_{{\cal C}}(S)$ be the minimum number of ${\cal C}$-edges of one-realizations of $S$. Then
  $$\delta_{{\cal C}}(S)=\left \{
\begin{array}{ll}2n_1-2n_s,  & \mbox{if~ $n_1-1,n_s+1\notin S$,} \\
2n_1-2n_s-2,  &\mbox{if~ $n_1-1,n_s+1\in S$,}\\
2n_1-2n_s-1,  &\mbox{otherwise.}\\
\end{array}
\right. $$
\end{thm}


\section{Proof of Theorem~\ref{thm1}}

 We first show that the number $\delta_{{\cal D}}(S)$ given in Theorem~\ref{thm1} is a lower bound on
the minimum number of ${\cal D}$-edges of one-realizations of $S$.

\begin{lemma}\label{lem2.1}
  $$\delta_{{\cal D}}(S)\geq\left \{
\begin{array}{ll}\frac{n_1(n_1-1)}{2}, & \mbox{if~ $n_1-1\notin S,$}\\
\frac{n_1(n_1-1)}{2}-1,&\mbox{if ~$n_1-1\in S.$}\\
\end{array}
\right. $$
 \end{lemma}

\proof Let ${\cal H}$ be a one-realization of $S$ and
$c=\{C_1,C_2,\ldots,C_{n_1}\}$ be a strict $n_1$-coloring of ${\cal
H}$.

\smallskip
 {\em Case 1.} $n_1-1\notin S$. If there exist two color classes $C_i$ and $C_j$
 such that ${\cal H}[C_i\cup C_j]$ has no ${\cal D}$-edges,
 then we can color the vertices in $C_i\cup C_j$ with a common color and get a strict $(n_1-1)$-coloring
  of ${\cal H}$, a contradiction. Hence, the first inequality holds.

\smallskip{\em Case 2.} $n_1-1\in S$. That is to say, $n_2=n_1-1$.
If there exist two distinct pairs of color classes, say $C_i,C_j$
and $C_k,C_l$, such that both ${\cal H}[C_i\cup C_j]$  and ${\cal
H}[C_k\cup C_l]$ have no ${\cal D}$-edges, then we can color the
vertices in $C_i\cup C_j$ with a common color and get a strict
$n_2$-coloring of ${\cal H}$; and color the vertices in $C_k\cup
C_l$ with a common color and get another strict $n_2$-coloring of
${\cal H}$, a contradiction. Hence, the second inequality
holds.$\qed$

Next we construct   mixed hypergraphs which meet  the bounds in
Lemma~\ref{lem2.1}. Let
$$
 \begin{array}{rcl}
{\cal D}_{n_1,\ldots,n_s}^{*}&=&\{\{(i,i,\ldots,i),(j,j,\ldots,j)\}|i\neq j, i,j\in [n_s-1]\}\\
&\cup&\{\{(i,\ldots,i),(x_1,\ldots,x_{s-1},n_s)\}|i\in [n_s-1]\}\\
&\cup&\{\{(x_1,\ldots,x_{s-1},1),(y_1,\ldots,y_{s-1},n_s)\}|x_i<y_i, i\in [s-1]\}.\\
\end{array}
$$
Then ${\cal H}_{n_1,\ldots,n_s}^*= (X_{n_1,\ldots,n_s},{\cal
C}_{n_1,\ldots,n_s}, {\cal D}_{n_1,\ldots,n_s}^*)$   is a spanning sub-hypergraph of  ${\cal H}_{n_1,\ldots,n_s}$ with $\frac{n_1(n_1-1)}{2}$ ${\cal D}$-edges.

By \cite[Theorem 2.4]{zdw1} all the strict colorings of ${\cal
H}_{n_1,\ldots,n_s}$ are as follows:
\begin{eqnarray}\label{eq1}
c_i^s=\{X_{i1}^s,X_{i2}^s,\ldots,X_{in_i}^s\}, \quad i\in [s],
\end{eqnarray}
 where $ X^s_{ij}=\{(x_1,x_2,\ldots,x_s)\in X_{n_1,\ldots,n_s}\mid x_i=j\},
j=1,2,\ldots,n_i.$
 Therefore, $c_1^s,\ldots,c_s^s$ are strict
colorings of ${\cal H}_{n_1,\ldots,n_s}^*$.

In the following we shall prove by induction on $s$  that
$c_1^s,\ldots,c_s^s$ are all the strict colorings of ${\cal
H}_{n_1,\ldots,n_s}^*$, which follows that  ${\cal
H}_{n_1,\ldots,n_s}^*$ is a one-realization of $S$.

\begin{lemma}\label{lem2.2}
 ${\cal H}_{n_1,n_2}^*$ is a one-realization of
$\{n_1, n_2\}$.
\end{lemma}

\proof Let $c=\{C_1,C_2,\ldots,C_m\}$ be a strict coloring of
${\cal H}_{n_1,n_2}^*$. The vertices $(1,1), (2,2), \ldots,
(n_2,n_2)$ fall into distinct color classes, say $(i,i)\in C_i, i\in [n_2]$. From the ${\cal C}$-edge $\{(n_2,1),(1,1),(n_2,n_2)\}$, we have
$(n_2,1)\in C_1\cup C_{n_2}$.

\smallskip
{\em Case 1.} $(n_2,1)\in C_1$. For any $k\in [n_1-n_2]$, from the ${\cal D}$-edge $\{(n_2,1),(n_2+k,n_2)\}$
and the ${\cal C}$-edge $\{(n_2+k,n_2),(n_2,1),(n_2,n_2)\}$, we have
$(n_2+k,n_2)\in C_{n_2}$. For any $k\in [n_1-n_2-1]$,  the ${\cal C}$-edge $
\{(n_2+k,1),(n_2,1),(n_2+k,n_2)\}$ and the ${\cal D}$-edge $\{(n_2+k,1),(n_1,n_2)\}$ imply that $(n_2+k,1)\in C_1$.
  Therefore, $c=c_2^2$.

\smallskip{\em Case 2.} $(n_2,1)\in C_{n_2}$.
The ${\cal D}$-edges $\{(n_2,1),(n_2+1,n_2)\},
\{(i,i),(n_2+1,n_2)\}, i\in [n_2-1]$ imply that $(n_2+1,n_2)\notin
C_i$ for any $i\in [n_2]$. Suppose $(n_2+1, n_2)\in C_{n_2+1}$.
Since  $\{(n_2+1,1),(n_2,1),(n_2+1,n_2)\},
\{(n_2+1,1),(n_2+1,n_2),(1,1)\}$ are ${\cal C}$-edges, $(n_2+1,1)\in
C_{n_2+1}$.
Similarly, $(n_2+k,n_2),(n_2+k,1)\in C_{n_2+k}$ for any $k\in
[n_1-n_2-1]$ and $(n_1,n_2)\in C_{n_1}$. Therefore, $c=c_1^2$.
$\qed$

\begin{thm}\label{thm2.3}
${\cal H}_{n_1,\ldots,n_s}^*$ is a one-realization of $S$.
\end{thm}

\proof By Lemma~\ref{lem2.2}   the conclusion is true for $s=2$.
Suppose that $s\geq 3$, and the conclusion is true for the case of
$s-1$. Let $X'=\{(x_1,x_2,\ldots,x_s)\in
X_{n_1,\ldots,n_s}|x_1=x_2\}$. Then ${\cal H}^*={\cal
H}_{n_1,\ldots,n_s}^*[X']$ is isomorphic to ${\cal
H}_{n_2,\ldots,n_s}^*$. By induction, all the strict colorings of
${\cal H}^*$ are as follows:
\begin{eqnarray}\label{eq2}
c'_i=\{X'_{i1},X'_{i2},\ldots, X'_{in_i}\}, \quad i\in [s]\setminus
\{1\},
\end{eqnarray}
 where $X'_{ij}=X'\cap X_{ij}^s, j\in [n_i]$.

For any strict coloring $c=\{C_1,C_2,\ldots,C_m\}$ of ${\cal
H}_{n_1,\ldots,n_s}^*$, the vertices $(1,1,\ldots,1),$
$(2,2,\ldots,2), \ldots, (n_s,n_s,\ldots,n_s)$  fall into distinct
color classes, say $(i,i,\ldots,i)\in C_i, i\in [n_s]$.  There are the
following two possible cases.

 \smallskip
{\em Case 1.}  $c|_{X'}=c_2'$. That is to say, for any $i\in [n_2]$ and any $(i,i,x_3,\ldots,x_s)\in X'$, $(i,i,x_3,\ldots,x_s)\in C_i$.
 The ${\cal C}$-edge $\{(1,1,\ldots,1),(n_2,1,\ldots,1),(n_2,n_2,n_3,\ldots,n_s)\}$
 implies that $(n_2,1,\ldots,1)\in C_1\cup C_{n_2}$.

\smallskip
 {\em Case 1.1.} $(n_2,1,\ldots,1)\in C_1$. For any $k\in [n_1-n_2]$, from the
${\cal C}$-edge $\{(n_2+k,n_2,n_3,\ldots,n_s),(n_2,1,\ldots,1),(n_2,n_2,n_3,\ldots,n_s)\}$
 and the ${\cal D}$-edge $\{(1,1,\ldots,1), (n_2+k,n_2,n_3,\ldots,n_s)\}$,
we observe $(n_2+k,n_2,n_3,\ldots,n_s)\in C_{n_2}$. For any $k\in
[n_1-n_2-1]$, by the ${\cal C}$-edge
$\{(n_2+k,1,\ldots,1),(n_2,n_2,n_3,\ldots,n_s),(n_2,1,\ldots,1)\}$
and the ${\cal D}$-edge
$\{(n_2+k,1,\ldots,1),(n_1,n_2,n_3,\ldots,n_s)\}$, we obtain
$(n_2+k,1,\ldots,1)\in C_1$.  Therefore, $c=c_2^s$.

\smallskip{\em Case 1.2.} $(n_2,1,\ldots,1)\in C_{n_2}$.
The fact that $\{(i,i,\ldots,i),(n_2+1,n_2,\ldots,n_s)\}$ is a
${\cal D}$-edge follows that $(n_2+1,n_2,\ldots,n_s)\notin C_i$ for
any $i\in [n_s-1]$. For each $t\in [s]\setminus \{1,2\}$ and $k\in
[n_{t-1}-n_t-1]\cup \{0\}$, the ${\cal D}$-edge
$\{(\underbrace{n_t+k,\ldots,n_t+k}_{t-1},1,\ldots,1),(n_2+1,n_2,\ldots,n_s)\}$
 implies that $(n_2+1,n_2,\ldots,n_s)\notin C_{n_t+k}$. Since $\{(n_2,1,\ldots,1),(n_2+1,n_2,\ldots,n_s)\}$
 is a ${\cal D}$-edge,  $(n_2+1,n_2,\ldots,n_s)\notin C_{n_2}$. Suppose $(n_2+1,n_2,\ldots,n_s)\in C_{n_2+1}$.
  From the ${\cal C}$-edges $\{(n_2+1,1,\ldots,1),(n_2+1,n_2,\ldots,n_s),(n_2,1,\ldots,1)\}$ and
   $\{(n_2+1,1,\ldots,1),(n_2+1,n_2,\ldots,n_s),(1,1,\ldots,1)\}$, we have $(n_2+1,1,\ldots,1)\in C_{n_2+1}$.
   Similarly, for any $k\in [n_1-n_2-1]$,  $(n_2+k,1,\ldots,1),(n_2+k,n_2,\ldots,n_s)\in
   C_{n_2+k}$
   and $(n_1,n_2,\ldots,n_s)\in C_{n_1}$. Hence,  $c=c_1^s$.

\smallskip
{\em Case 2.}  $c|_{X'}=c_t'$ for some   $t\in [s]\setminus
\{1,2\}$.  In this case, we observe $(n_2-1,n_2-1,n_3,\ldots,n_s)\in
C_{n_t}$ and $(n_2-1,n_2-1,1,\ldots,1)\in C_1$. For any $k\in
[n_1-n_2]$, from the ${\cal C}$-edge
$\{(n_2+k,n_2,\ldots,n_s),(n_2-1,n_2-1,n_3,\ldots,n_s),(n_2-1,n_2-1,1,\ldots,1)\}$
and the ${\cal D}$-edge $\{(1,1,\ldots,1),(n_2+k,n_2,\ldots,n_s)\}$,
we have $(n_2+k,n_2,\ldots,n_s)\in C_{n_t}$. Furthermore,  the
${\cal C}$-edge
$\{(n_2+k,1,\ldots,1),(n_2+k,n_2,\ldots,n_s),(1,1,\ldots,1)\}$ and
the ${\cal D}$-edge $\{(n_2+k,1,\ldots,1),(n_1,n_2,\ldots,n_s)\}$
imply that $(n_2+k,1,\ldots,1)\in C_1$ for any $k\in [n_1-n_2-1]\cup
\{0\}$.   Therefore, $c=c_t^s$.$\qed$

When $n_1-1\in S$, that is to say, $n_2=n_1-1$, let
$X''=X_{n_1,\ldots,n_s}\setminus \{(n_2,1,\ldots,1)\}$.   Then
 ${\cal G}_{n_1,\ldots,n_s}^*={\cal H}_{n_1,\ldots,n_s}^*[X'']$ is a
spanning sub-hypergraph of ${\cal G}_{n_1,\ldots,n_s}$ with
$\frac{n_1(n_1-1)}{2}-1$ ${\cal D}$-edges. Note that for any $i\in
[s]$,
\begin{eqnarray}\label{eq3}
c_i''=\{X_{i1}'',X_{i2}'',\ldots,X_{in_i}''\}
\end{eqnarray}
 is a strict
$n_i$-coloring of ${\cal G}_{n_1,\ldots,n_s}^*$, where
$X_{ij}''=X''\cap X_{ij}^s, j=1,2,\ldots,n_i$.

 \begin{thm}\label{thm2.4} If $n_1-1\in S$, then ${\cal G}_{n_1,\ldots,n_s}^*$ is a one-realization of $S$.
\end{thm}

\proof Referring  to the proof of Theorem~\ref{thm2.3},  all the
strict colorings of ${\cal H}^*$ are $c_2',c_3',\ldots,c_s'$. For
any strict coloring $c=\{C_1,C_2,\ldots,C_m\}$ of ${\cal
G}_{n_1,\ldots,n_s}^*$, we focus on the restriction of $c$ on $X'$
and get the following two possible cases.

\smallskip{\em Case 1.}  $c|_{X'}=c_2'$.
That is to say, $(i,i,x_3,\ldots,x_s)\in C_i$ for any $i\in [n_2]$
and $(i,i,x_3,\ldots,x_s)\in X'$.   By the proof of
Theorem~\ref{thm2.3}  we have $(n_1,n_2,n_3,\ldots,n_s)\notin C_j$
for any $j\in [n_2-1]$. If
 $(n_1,n_2,n_3,\ldots,n_s)\in C_{n_2},$ then $c=c_2''$; If  $(n_1,n_2,n_3,\ldots,n_s)\notin C_{n_2},$ then
$(n_1,n_2,n_3,\ldots,n_s)\in C_{n_1}$, which implies  that
 $c=c_1''$.

\smallskip{\em Case 2.}  $c|_{X'}=c_t'$ for some $t\in [s]\setminus \{1,2\}$. That is to say,
$(n_t,\ldots,n_t,n_{t+1},\ldots,n_s)\in C_{n_t}$ and $(\underbrace{n_t,\ldots,n_t}_{t-1},1,\ldots,1)\in C_1$. The ${\cal D}$-edge $\{(1,1,\ldots,1),(n_1,n_2,\ldots,n_s)\}$ and the  ${\cal C}$-edge $\{(n_1,n_2,\ldots,n_s),(n_t,\ldots,n_t,n_{t+1},\ldots,n_s),(n_t,\ldots,n_t,1,\ldots,1)\}$ imply that $(n_1,n_2,\ldots,n_s)\in C_{n_t}$. Hence, $c=c_t''$.

By above discussion, $c_1'', c_2'',\ldots,c_s''$ are all the strict
colorings of ${\cal G}_{n_1,\ldots,n_s}^*$. Therefore, ${\cal
G}_{n_1,\ldots,n_s}^*$ is a one-realization of $S$. $\qed$

 Combining Lemma~\ref{lem2.1}, Theorems~\ref{thm2.3} and  \ref{thm2.4}, we complete the proof of
 Theorem~\ref{thm1}.

\section{Proof of Theorem~\ref{thm2}}

 In \cite{zdw1}, we determined the minimum number of vertices of
one-realizations of $S$.

\begin{thm}\label{thm3.0}{\rm(\cite{zdw1})} Let $\delta(S)$ denote the minimum number of vertices of one-realizations of $S$. Then
  $$\delta(S)=\left \{
\begin{array}{ll}2n_1-n_s,  & \mbox{if~ $n_1-1\notin S,$} \\
2n_1-n_s-1,  &\mbox{if~ $ n_1-1\in S.$}\\
\end{array}
\right. $$
\end{thm}

 Let  ${\cal H}$ be a mixed hypergraph with a strict $m$-coloring
$c=\{C_1,C_2,\ldots, C_m\}$.
 For any $i$,   the {\em pair graph}  $G_i$ of $C_i$ under $c$ has $C_i$ as its vertex set,
 and two vertices $x$ and $y$ are adjacent if there exists a $\cal
 C$-edge $C$ such that   $C\cap C_i= \{x,y\}$ and  $|C\cap C_j|\leq 1$
 for any $j\in [m]\setminus \{i\}$. Note that the definition of pair
 graphs  is different from that in \cite{Diao1}.

The following result shows that the number $\delta_{{\cal C}}(S)$
given in Theorem 1.2 is a lower bound on the minimum number of
${\cal C}$-edges of one-realizations of $S$.

 \begin{lemma}\label{lem3.1}
  $$\delta_{{\cal C}}(S)\geq\left \{
\begin{array}{ll}\delta(S)-n_s,  & \mbox{if~ $n_s+1\notin S,$} \\
\delta(S)-n_s-1,  &\mbox{if~ $ n_s+1\in S.$}\\
\end{array}
\right. $$
 \end{lemma}

\proof Let ${\cal H}=(X, {\cal C}, {\cal D})$ be  a one-realization
of $S$ and  $c=\{C_1,C_2,\ldots,
C_{n_s}\}$ be a strict $n_s$-coloring of ${\cal H}$. For any $i\in
[n_s]$, let $G_i=(C_i, E_i)$ be the pair graph of $C_i$ under $c$.

\smallskip
 {\em Case 1.} $n_s+1\notin S$. If   $G_i$ is a disconnected graph with $G_i^1=(C_i^1, E_i^1)$
 as one of its connected components, then $\{C_1,\ldots,C_{i-1},C_i^1, (C_i\setminus C_i^1),C_{i+1},\ldots,C_{n_s}\}$
  is a strict $(n_s+1)$-coloring of ${\cal H}$, a contradiction. Therefore, each   $G_i$ is connected; and so $|E_i|\geq |C_i|-1$.
 Since every ${\cal C}$-edge of ${\cal H}$ corresponds to at most one edge in the $n_s$ pair graphs,
  we have $\delta_{{\cal C}}(S)\geq (|C_1|-1)+(|C_2|-1)+\cdots+ (|C_{n_s}|-1)= |X|-n_s\geq \delta(S)-n_s$.

\smallskip
 {\em Case 2.} $n_s+1\in S$.
Assume that there exist two  disconnected pair graphs, say
  $G_1$ and $G_2$. Pick     a connected component $G_1^1=(C_1^1, E_1^1)$ of $G_1$,
   and   a connected component $G_2^1=(C_2^1,E_2^1)$ of $G_2$.
   Then $\{C_1^1, (C_1\setminus C_1^1),C_2,\ldots, C_{n_s}\}$ and $\{C_1,C_2^1,
   (C_2\setminus C_2^1),C_3,\ldots, C_{n_s}\}$ are two distinct strict
   $(n_s+1)$-colorings of ${\cal H}$,
   a contradiction.   If there exists a   pair graph, say $G_1$, with at
   least three connected components  $G_1^1=(C_1^1, E_1^1)$, $G_1^2=(C_1^2, E_1^2)$
   and $G_1^3=(C_1^3, E_1^3)$, then $\{C_1^1, (C_1\setminus C_1^1),C_2,\ldots,C_{n_s}\}$,
     $\{C_1^2, (C_1\setminus C_1^2),C_2,\ldots,C_{n_s}\}$
     are two distinct strict $(n_s+1)$-colorings of ${\cal H}$, a contradiction.
     Hence, there exists at most one   pair graph with   two connected components.
       It follows that ${\cal H}$ has at least
       $\delta(S)-n_s-1$ ${\cal C}$-edges. $\qed$

Next we   construct   mixed hypergraphs which meet  the bounds in
Lemma~\ref{lem3.1}. Let
$$\begin{array}{l}
{\cal
C}_{n_1,\ldots,n_s}^{{\star}}\\
=\bigcup\limits_{t=2}^s\bigcup\limits_{j=n_t+1}^{n_{t-1}-1}\{\{(j,\ldots,j,n_t,\ldots,n_s),(\underbrace{j,\ldots,j}_{t-1},1,\ldots,1),
(\underbrace{j-1,\ldots,j-1}_{t-1},1,\ldots,1)\}\}\\
\cup\bigcup\limits_{t=2}^s\bigcup\limits_{j=n_t}^{n_{t-1}-1}\{(\{\underbrace{j,\ldots,j}_{t-1},1,\ldots,1),
(j,\ldots,j,n_t,\ldots,n_s),(j+1,\ldots,j+1,n_t,\ldots,n_s)\}\}\\
\cup\bigcup\limits_{t=2}^s\{\{(n_t,\ldots,n_t,n_t,\ldots,n_s),(\underbrace{n_t,\ldots,n_t}_{t-1},1,\ldots,1),
(1,1,\ldots,1)\}\}.
\end{array}
$$
Then ${\cal H}_{n_1,\ldots,n_s}^{{\star}}= (X_{n_1,\ldots,n_s},{\cal
C}_{n_1,\ldots,n_s}^{{\star}}, {\cal D}_{n_1,\ldots,n_s})$   is a
spanning sub-hypergraph of ${\cal H}_{n_1,\ldots,n_s}$ with
$2n_1-2n_s$ ${\cal C}$-edges. Therefore, $c_1^s, c_2^s, \ldots,
c_s^s$ in (\ref{eq1}) are strict colorings of ${\cal
H}_{n_1,\ldots,n_s}^{\star}$.

In the following, we shall prove by induction on $s$  that
$c_1^s,\ldots,c_s^s$ are all the strict colorings of ${\cal H}_{n_1,\ldots,n_s}^{\star}$, which follows that ${\cal H}_{n_1,\ldots,n_s}^{\star}$ is a one-realization of $S$.

\begin{lemma}\label{lem3.2}  ${\cal H}_{n_1,n_2}^{\star}$ is a one-realization of
$\{n_1, n_2\}$.
\end{lemma}

\proof Suppose $c=\{C_1,C_2,\ldots,C_m\}$ is a strict coloring of
${\cal H}_{n_1,n_2}^{\star}$. The vertices $(1,1), (2,2), \ldots,
(n_2,n_2)$ fall into distinct color classes, say $(i,i)\in C_i, i\in
[n_2]$. From the ${\cal C}$-edge $\{(n_2,1),(n_2,n_2),(1,1)\}$, we
have $(n_2,1)\in C_1\cup C_{n_2}$.

\smallskip
{\em Case 1.} $(n_2,1)\in C_1$.  The ${\cal C}$-edge
$\{(n_2,1),(n_2,n_2),(n_2+1,n_2)\}$ and the ${\cal D}$-edge
$\{(1,1),(n_2+1,n_2)\}$  imply that $(n_2+1,n_2)\in C_{n_2}$. From
the ${\cal D}$-edge $\{(n_2+1,1),(n_2,n_2)\}$ and the ${\cal
C}$-edge $\{(n_2+1,1),(n_2+1,n_2),(n_2,1)\}$, we have $(n_2+1,1)\in
C_1$. Similarly, for any $k\in [n_1-n_2-1]$, we have $(n_2+k,1)\in
C_1, (n_2+k,n_2)\in C_{n_2}$
 and $(n_1,n_2)\in C_{n_2}$. Therefore, $c=c_2^2$.

\smallskip{\em Case 2.} $(n_2,1)\in C_{n_2}$. Observe $(n_2+1,n_2)\notin C_i$ for each $i\in [n_2-1]$.
The ${\cal D}$-edge $\{(n_2,1),(n_2+1,n_2)\}$ implies that
$(n_2+1,n_2)\notin C_{n_2}$. Suppose $(n_2+1, n_2)\in C_{n_2+1}$.
From the ${\cal C}$-edge $\{(n_2+1,1),(n_2+1,n_2),(n_2,1)\}$ and the
${\cal D}$-edge $\{(n_2+1,1),(n_2,n_2)\}$, we get  $(n_2+1,1)\in
C_{n_2+1}$. Similarly, $(n_2+k,n_2),(n_2+k,1)\in C_{n_2+k}$ for any
$k\in [n_1-n_2-1]$ and $(n_1,n_2)\in C_{n_1}$. Therefore, $c=c_1^2$.
$\qed$

\begin{thm}\label{thm3.3}
${\cal H}_{n_1,\ldots,n_s}^{\star}$ is a one-realization of $S$.
\end{thm}

\proof By Lemma~\ref{lem3.2}   the conclusion is true for $s=2$.
Suppose $s\geq 3$, and the conclusion is true for the case of $s-1$.
Let $X'=\{(x_1,x_2,,\ldots,x_s)\in X_{n_1,\ldots,n_s}|x_1=x_2\}$.
Then ${\cal H}^{\star}={\cal H}_{n_1,\ldots,n_s}^{\star}[X']$ is
isomorphic to ${\cal H}_{n_2,n_3,\ldots,n_s}^{\star}$. By induction,
all the strict colorings of ${\cal H}^{\star}$ are $c'_2,
c_3',\ldots,c_s'$ in (\ref{eq2}).

For any strict coloring $c=\{C_1,C_2,\ldots,C_m\}$ of ${\cal
H}_{n_1,\ldots,n_s}^{\star}$, the vertices $(1,1,\ldots,1),$
$(2,2,\ldots,2), \ldots, (n_s,n_s,\ldots,n_s)$  fall into distinct
color classes, say $(i,i,\ldots,i)\in C_i, i\in [n_s]$.  We focus on the restriction of $c$ on $X'$ and get the
following two possible cases.

 \smallskip
{\em Case 1.}  $c|_{X'}=c_2'$. That is to say, for any $i\in [n_2]$ and any $(i,i,x_3,\ldots,x_s)\in X'$, $(i,i,x_3,\ldots,x_s)\in C_i$.
 The ${\cal C}$-edge $\{(n_2,1,\ldots,1),(n_2,n_2,n_3,\ldots,n_s),(1,1,\ldots,1)\}$ implies that $(n_2,1,\ldots,1)\in C_1\cup C_{n_2}$.

\smallskip
 {\em Case 1.1.} $(n_2,1,\ldots,1)\in C_1$. From the ${\cal D}$-edge $\{(1,1,\ldots,1), (n_2+1,n_2,n_3,\ldots,n_s)\}$  and the
${\cal C}$-edge $\{(n_2,1,\ldots,1),(n_2,n_2,n_3,\ldots,n_s),(n_2+1,n_2,n_3,\ldots,n_s)\}$,
we observe $(n_2+1,n_2,n_3,\ldots,n_s)\in C_{n_2}$. Then the ${\cal D}$-edge $\{(n_2+1,1,\ldots,1),(n_2,n_2,n_3,\ldots,n_s)\}$ and the ${\cal C}$-edge $\{(n_2+1,1,\ldots,1),(n_2+1,n_2,n_3,\ldots,n_s),(n_2,1,\ldots,1)\}$  imply that $(n_2+1,1,\ldots,1)\in C_1$.
Similarly, for any $k\in [n_1-n_2-1]$,
$(n_2+k,n_2,\ldots,n_s)\in C_{n_2},(n_2+k,1,\ldots,1)\in C_1$
and $(n_1,n_2,n_3,\ldots,n_s)\in C_{n_2}$.  Therefore, $c=c_2^s$.

\smallskip{\em Case 1.2.} $(n_2,1,\ldots,1)\in C_{n_2}$.
The fact that $\{(n_2+1,n_2,\ldots,n_s),(i,i,\ldots,i)\}$ is a
${\cal D}$-edge follows that $(n_2+1,n_2,\ldots,n_s)\notin C_i$ for
any $i\in [n_s-1]$. For any $t\in [s]\setminus \{1,2\}$ and $k\in
[n_{t-1}-n_t-1]\cup \{0\}$, the ${\cal D}$-edge
$\{(n_2+1,n_2,\ldots,n_s),(\underbrace{n_t+k,\ldots,n_t+k}_{t-1},1,\ldots,1)\}$
implies that $(n_2+1,n_2,\ldots,n_s)\notin C_{n_t+k}$. Since
$\{(n_2+1,n_2,\ldots,n_s),(n_2,1,\ldots,1)\}$ is a ${\cal D}$-edge,
$(n_2+1,n_2,\ldots,n_s)\notin C_{n_2}$. Suppose
$(n_2+1,n_2,\ldots,n_s)\in C_{n_2+1}$. From the ${\cal C}$-edge
$\{(n_2+1,1,\ldots,1),(n_2+1,n_2,\ldots,n_s),(n_2,1,\ldots,1)\}$ and
the ${\cal D}$-edge $\{(n_2+1,1,\ldots,1),(n_2,n_2,\ldots,n_s)\}$,
we have $(n_2+1,1,\ldots,1)\in C_{n_2+1}$. Similarly,
$(n_2+k,1,\ldots,1),(n_2+k,n_2,\ldots,n_s)\in C_{n_2+k}$ for any
$k\in [n_1-n_2-1]$ and  $(n_1,n_2,\ldots,n_s)\in C_{n_1}$.
Therefore,  $c=c_1^s$.

\smallskip
  {\em Case 2.}  $c|_{X'}=c_t'$ for some $t\in [s]\setminus
\{1,2\}$. Then    $(n_2,n_2,n_3,\ldots,n_s),
(n_2-1,n_2-1,n_3,\ldots,n_s)\in C_{n_t},(n_2-1,n_2-1,1,\ldots,1)\in
C_1$. From the ${\cal C}$-edge $\{(n_2,1,\ldots,1),$
$(n_2,n_2,n_3,\ldots,n_s),(1,1,\ldots,1)\}$ and the ${\cal D}$-edge
$\{(n_2,1,\ldots,1), (n_2-1,n_2-1,n_3,\ldots,n_s)\}$, we have
$(n_2,1,\ldots,1)\in C_1$. From the  ${\cal C}$-edge
$\{(n_2,1,\ldots,1),(n_2,n_2,n_3,\ldots,n_s),(n_2+1,n_2,\ldots,n_s)\}$
and the ${\cal D}$-edge $\{(n_2+1,n_2,\ldots,n_s),(1,1,\ldots,1)\}$,
we obtain $(n_2+1,n_2,\ldots,n_s)\in C_{n_t}$. Then the ${\cal
C}$-edge
$\{(n_2+1,1,\ldots,1),(n_2+1,n_2,\ldots,n_s),(n_2,1,\ldots,1)\}$ and
the ${\cal D}$-edge $\{(n_2,n_2,\ldots,n_s),(n_2+1,1,\ldots,1)\}$
imply that $(n_2+1,1,\ldots,1)\in C_1$. Similarly, for any $k\in
[n_1-n_2-1]$, $(n_2+k,n_2,\ldots,n_s)\in C_{n_t},
(n_2+k,1,\ldots,1)\in C_1$, and $(n_1,n_2,\ldots,n_s)\in C_{n_t}$.
Hence, $c=c_t^s$. $\qed$

When $n_1-1\in S$, that is to say, $n_1=n_2+1$, let ${\cal
G}_{n_1,\ldots,n_s}^{\star}={\cal
H}_{n_1,\ldots,n_s}^{\star}[X'']+e_1,$  where
$e_1=\{(n_2-1,n_2-1,1,\ldots,1),(n_2-1,n_2-1,n_3,\ldots,n_s),(n_1,n_2,\ldots,n_s)\}$
and $X''=X_{n_1,\ldots,n_s}\setminus \{(n_2,1,\ldots,1)\}.$ Then,
${\cal G}_{n_1,\ldots,n_s}^{\star}$ is a spanning sub-hypergraph of
${\cal G}_{n_1,\ldots,n_s}$ with $2n_1-2n_s-1$ ${\cal C}$-edges.
Therefore, $c_1'', c_2'', \ldots, c_s''$ in (\ref{eq3}) are strict
colorings of ${\cal G}_{n_1,\ldots,n_s}^{\star}$.

 \begin{thm}\label{thm3.4} If $n_1=n_2+1$, then ${\cal G}_{n_1,\ldots,n_s}^{\star}$ is a one-realization of $S$.
\end{thm}

\proof Let  $X'=\{(x_1,x_2,\ldots,x_s)\in
X_{n_1,\ldots,n_s}|x_1=x_2\}$ and ${\cal G}^{\star}={\cal
G}_{n_1,\ldots,n_s}^{\star}[X']$. Referring  to the proof of
Theorem~\ref{thm3.3},  all the strict colorings of ${\cal
H}^{\star}$ are $c'_2, c_3',\ldots,c_s'$ in (\ref{eq2}). For any
strict coloring $c=\{C_1,C_2,\ldots,C_m\}$ of ${\cal
G}_{n_1,\ldots,n_s}^{\star}$, we focus on the restriction of $c$ on
$X'$ and get the following two possible cases.

\smallskip{\em Case 1.}  $c|_{X'}=c_2'$.
That is to say, $(i,i,x_3,\ldots,x_s)\in C_i$ for any $i\in [n_2]$ and $(i,i,x_3,\ldots,x_s)\in X'$.
 By the same discussion as we do in Theorem~\ref{thm3.3}, we have $(n_1,n_2,n_3,\ldots,n_s)\notin C_j$ for any $j\in [n_2-1]$.
  If $(n_1,n_2,n_3,\ldots,n_s)\in C_{n_2},$ then
  $c=c_2''$; If   $(n_1,n_2,n_3,\ldots,n_s)\notin C_{n_2},$ then
$(n_1,n_2,n_3,\ldots,n_s)\in C_{n_1}$ which follows that $c=c_1''$.

\smallskip   {\em Case 2.} $c|_{X'}=c_t'$ for some $t\in [s]\setminus \{1,2\}$.
Then $(n_2-1,n_2-1,n_3,\ldots,n_s)\in C_{n_t}$ and
$(1,1,\ldots,1),(n_2-1,n_2-1,1,\ldots,1)\in C_1$. From the ${\cal
C}$-edge
$\{(n_2-1,n_2-1,1,\ldots,1),(n_2-1,n_2-1,n_3,\ldots,n_s),(n_1,n_2,\ldots,n_s)\}$
and the ${\cal D}$-edge $\{(n_1,n_2,\ldots,n_s),$ $
(1,1,\ldots,1)\}$, we observe $(n_1,n_2,\ldots,n_s)\in C_{n_t}$.
Hence, $c=c_t''$. $\qed$

When $n_1-1\notin S$ and $n_s+1\in S$,   let
 ${\cal H}_{n_1,\ldots,n_s}^{{\star}1}={\cal
 H}_{n_1,\ldots,n_{s-2},n_s}^{\star}-e_2$,
 where $e_2=\{(n_s,\ldots,n_s,1),(n_s,\ldots,n_s,n_s),(1,\ldots,1,1)\}.$
Then ${\cal H}_{n_1,\ldots,n_s}^{{\star}1}$ is a spanning
sub-hypergraph of ${\cal H}_{n_1,\ldots,n_{s-2},n_s}^{\star}$ with
$2n_1-2n_s-1$ ${\cal C}$-edges. Note that, for any $i\in
[s]\setminus \{s-1\}$,
$$c_i^{s-1}=\{X_{i1}^{s-1},X_{i2}^{s-1},\ldots,X_{in_i}^{s-1}\}$$ is
a strict $n_i$-coloring of ${\cal H}_{n_1,\ldots,n_s}^{{\star}1}$,
where $X^{s-1}_{ij}=\{(x_1,\ldots,x_{s-2},x_s)\in
X_{n_1,\ldots,n_{s-2},n_s}|~x_i=j\}, j=1,2,\ldots,n_i$.

   \begin{thm}\label{thm3.5} If $n_1-1\notin S$ and $n_s+1\in S$, then ${\cal H}_{n_1,\ldots,n_s}^{{\star}1}$ is a one-realization of $S$.
\end{thm}

\proof Let $c=\{C_1,C_2,\ldots,C_m\}$ be a strict coloring of ${\cal
H}_{n_1,\ldots,n_s}^{{\star}1}$. Then the vertices
$(i,i,\ldots,i),i=1,2,\ldots,n_s$ fall into distinct color classes.
Suppose $(i,i,\ldots,i)\in C_i$ for any $i\in [n_s]$. Note that
$(n_s,\ldots,n_s,1)\notin C_i, i\in [n_s]\setminus \{1, n_s\}$, we
get the following two possible cases.

 \smallskip
{\em Case 1.} $(n_s,\ldots,n_s,1)\in C_1\cup  C_{n_s}$. By the same
discussion  in Theorem~\ref{thm3.3}, we have $c\in \{c_1^{s-1},\ldots,c_{s-2}^{s-1},c_s^{s-1}\}$.

\smallskip
{\em Case 2.} $(n_s,\ldots,n_s,1)\notin C_i,i\in [n_s]$. Suppose
$(n_s,\ldots,n_s,1)\in C_{n_{s-1}}$. Then the ${\cal D}$-edge
$\{(n_s,\ldots,n_s,1),(n_s+1,\ldots,n_s+1,n_s)\}$ and the ${\cal
C}$-edge
$\{(n_s+1,\ldots,n_s+1,n_s),(n_s,\ldots,n_s,1),(n_s,\ldots,n_s,n_s)\}$
imply that $(n_s+1,\ldots,n_s+1,n_s)\in C_{n_s}$. From the ${\cal
C}$-edge
$\{(n_s+1,\ldots,n_s+1,1),(n_s+1,\ldots,n_s+1,n_s),(n_s,\ldots,n_s,1)\}$
and the ${\cal D}$-edge
$\{(n_s,\ldots,n_s,n_s),(n_s+1,\ldots,n_s+1,1)\}$, we have
$(n_s+1,\ldots,n_s+1,1)\in C_{n_{s-1}}$. Similarly, for any $k\in
[n_{s-2}-n_s-1]$,  $(n_s+k,\ldots,n_s+k,n_s)\in C_{n_s},
(n_s+k,\ldots,n_s+k,1)\in C_{n_{s-1}}$ and
$(n_{s-2},\ldots,n_{s-2},n_s)\in C_{n_s}$. Since
$\{(n_{s-2},\ldots,n_{s-2},1,1),(n_s,\ldots,n_s,n_s)\}$ is a ${\cal
D}$-edge and
$\{(n_{s-2},\ldots,n_{s-2},1,1),(n_{s-2},\ldots,n_{s-2},n_s),(1,\ldots,1)\}$
is a ${\cal C}$-edge, we get $(n_{s-2},\ldots,n_{s-2},1,1)\in C_1$.
Similar to the discussion   in Theorem~\ref{thm3.3}, we have
$(x_1,\ldots,x_{s-3},1,1)\in C_1$ and
$(x_1,\ldots,x_{s-3},n_{s-2},n_s)\in C_{n_s}$ for any
$(x_1,\ldots,x_{s-3},1,1)$, $(x_1,\ldots,x_{s-3},n_{s-2},n_s)\in
X_{n_1,\ldots,n_{s-2},n_s}$. Therefore, $c$ is a strict
$n_{s-1}$-coloring of ${\cal H}_{n_1,\ldots,n_s}^{{\star}1}$.

By the above discussion, we have that ${\cal
H}_{n_1,\ldots,n_s}^{{\star}1}$ is a one-realization of $S$.$\qed$

When $n_1-1, n_s+1\in S$, let
$X'''=X_{n_1,\ldots,n_{s-2},n_s}\setminus \{(n_2,1,\ldots,1)\}$ and
${\cal H}_{n_1,\ldots,n_s}^{{\star}2}={\cal
H}_{n_1,\ldots,n_s}^{{\star}1}[X''']+e_3$, where $e_3=
\{(n_2-1,n_2-1,1,\ldots,1),(n_2-1,n_2-1,n_3,\ldots,n_{s-2},n_s)$,
$(n_1,\ldots,n_{s-2},n_s)\}.$ By the discussion    in
Theorems~\ref{thm3.3}-\ref{thm3.5}, we have the following result.

  \begin{thm}\label{thm3.6} If $n_1-1, n_s+1\in S$, then ${\cal H}_{n_1,\ldots,n_s}^{{\star}2}$ is a one-realization of $S$.
\end{thm}

Combining  Theorem~\ref{thm3.0},
Lemma~\ref{lem3.1}, Theorems~\ref{thm3.3}-\ref{thm3.6}, we complete
the proof of Theorem~\ref{thm2}.

\section*{Acknowledgment}

 The research
 is supported by NSF of
Shandong Province (No. ZR2009AM013), NCET-08-0052, NSF of China
(10871027) and the Fundamental Research Funds for the Central
Universities of China.

\end{CJK*}

\end{document}